\documentclass[12pt]{article} 

\usepackage[utf8]{inputenc} 

\usepackage{geometry} 
\usepackage{graphicx} 

\usepackage{booktabs} 
\usepackage{array} 
\usepackage{paralist} 
\usepackage{verbatim} 
\usepackage{subfig} 

\usepackage{fancyhdr} 
\usepackage{sectsty}
\allsectionsfont{\sffamily\mdseries\upshape} 

\usepackage[nottoc,notlof,notlot]{tocbibind} 
\usepackage[titles,subfigure]{tocloft} 

\usepackage{url} 
\usepackage{amsthm, amssymb,amsmath} 

\title{Induced Subgraphs of Order Seven and Their Frequencies in $srg(n,k,1,2)$}
\author{Reimbay Reimbayev}
\date{} 

\begin{document}
\maketitle

\begin{abstract}
In this paper, we examine the structure of strongly regular graphs with parameters $\lambda = 1$ and $\mu = 2$. In particular, we provide a complete classification of induced subgraphs of order seven and determine their relative frequencies. These findings contribute to a finer understanding of the local structure of such graphs and may be useful in related combinatorial and algebraic investigations.
\end{abstract}

\section{Introduction}

Strongly regular graphs form an important class of combinatorial structures, arising naturally in algebraic graph theory, finite geometry, and design theory. They are characterized by four parameters $(n,k,\lambda,\mu)$, which impose strict regularity conditions on the adjacency relations between vertices. A graph with parameters $(n,k,\lambda,\mu)$ is said to be strongly regular if it is $k$-regular of order $n$, such that every pair of adjacent vertices has exactly $\lambda$ common neighbors, and every pair of non-adjacent vertices has exactly $\mu$ common neighbors \cite{Gordon, Brouwer}. Another definition of strongly regular graphs, perhaps more precise as it cuts away some trivial cases such as complete graphs, is given by means of spectral graph theory, in which the finite graph is strongly regular if its spectrum consists of exactly three eigenvalues, one of which is $k$ with multiplicity one \cite{BrouwerMaldeghem}.

Since their introduction by R.C. Bose in their foundational paper \cite{Bose}, the strongly regular graphs have been a central object of study in algebraic graph theory. Bose established many of the basic parameter relations for strongly regular graphs and highlighted their deep connections with combinatorial design theory. A systematic treatment of strongly regular graphs and their algebraic properties was later developed in the monograph Algebraic Graph Theory by Chris Godsil and Gordon Royle \cite{Godsil}, which remains a standard reference in the field. Earlier, spectral techniques were extensively developed by Dragoš Cvetković at al in Spectra of Graphs \cite{Cvetcovic}, providing powerful tools for analyzing eigenvalue constraints and non-existence results.

A large body of research has focused on existence, classification, and construction of strongly regular graphs with prescribed parameters. Surveys and databases compiled by Andries Brouwer and collaborators remain essential tools in this direction \cite{Brouwer}.

The specific case for $\lambda = 1$ and $\mu = 2$ has attracted attention due to its connection with long-standing existence problems. For instance, the possible existence of a graph with parameters $\mathrm{srg}(99,14,1,2)$ has remained open for decades \cite{Conway}. But this is just one graph from the family of strongly regular graphs for which only few known to exist, e.g. those with valencies $k=2, 4$ and, surprisingly 22 \cite{Berlekamp}. There have been some extensive studies by Makhnev et. al. on the structure of such graphs with regard to their automorphism groups \cite{MakhnevMinkova}. Using Wilbrink and Brouwer's lemma \cite{WilbrinkBrouwer}, Lou and Murin were able to establish a forbidden subgraph of order 9 in case when $k=14$ \cite{LouMurin}.

Our recent work investigates structural properties of this family of graphs and establishes lower bounds on the number of hexagons, linking cycle structure to existence questions \cite{ReiLowerBound}. This demonstrates how local subgraph counts can provide insight into global feasibility.
In a similar vein, we have studied induced subgraphs of order six in strongly regular graphs with $\lambda = 1$ and $\mu = 2$, giving a complete enumeration of all such configurations \cite{ReiAllSix}. Furthermore, we have also provided complete enumeration of induced Hamiltonian subgraphs of order seven with their relative frequencies \cite{ReiHamiltonianSeven}. This line of research emphasizes the importance of understanding small subgraph structure as a means of probing the fine combinatorial properties of these graphs. 

The present work builds on these developments by extending the local analysis to subgraphs of order seven and determining their relative frequencies, thereby contributing to a more detailed understanding of the structure of this important class of strongly regular graphs.

\section{Subgraphs and Frequencies}

\begin{figure}
	\includegraphics[width=1.0\textwidth]{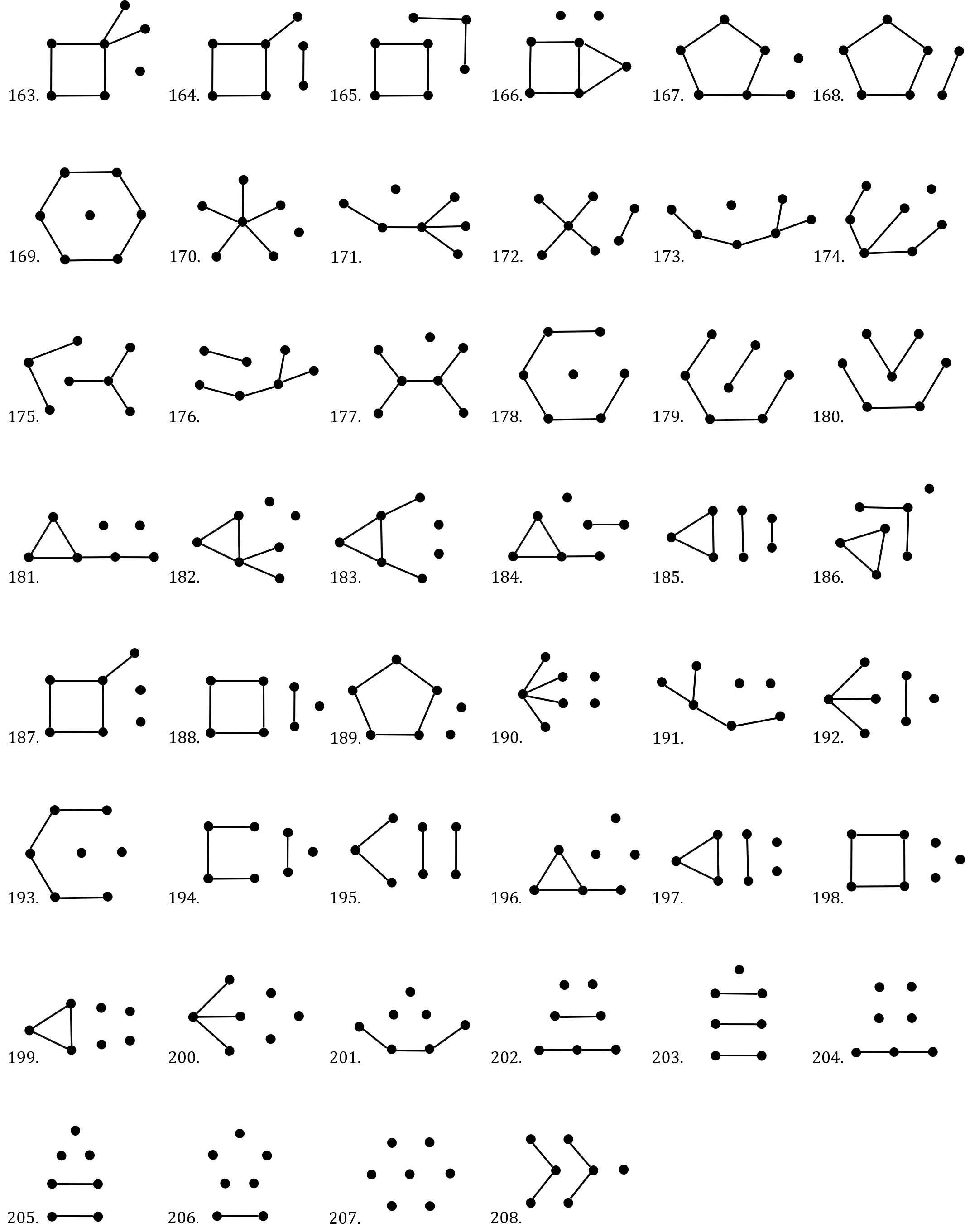}
		\centering
		\caption{All possible subgraphs of order six in $srg(n,k,1,2)$.}
		\label{mainFigure}
\end{figure}

There are 208 subgraphs of order seven in $srg(n,k,1,2)$. All of them are given in the main figure - Figure \ref{mainFigure}, spanned over several pages below. Their frequencies are given further, denoted by $z_i$, while the subgraphs are denoted by capital letter $Z_i$ . Some of the subgraphs might be non existent, i.e. have a frequency zero depending on the values of $n_3$ and $z_{11}$. Theses are the free variables as the exact values for $z_i$-s are not known. We follow the same convention as in previous work \cite{ReiAllSix} for $n_3$, while $z_{11}$ is the frequency of the graph $Z_{11}$ as agreed. We omit all the derivations of the formulas for frequencies as they are very similar to what have done earlier \cite{ReiHamiltonianSeven}. In fact, this is the extension on the latter work, in which we have looked at Hamiltonian subgraphs of order seven, while here we have considered all the rest of the subgraphs and their frequencies.

\begin{figure}
	\includegraphics[width=1.0\textwidth]{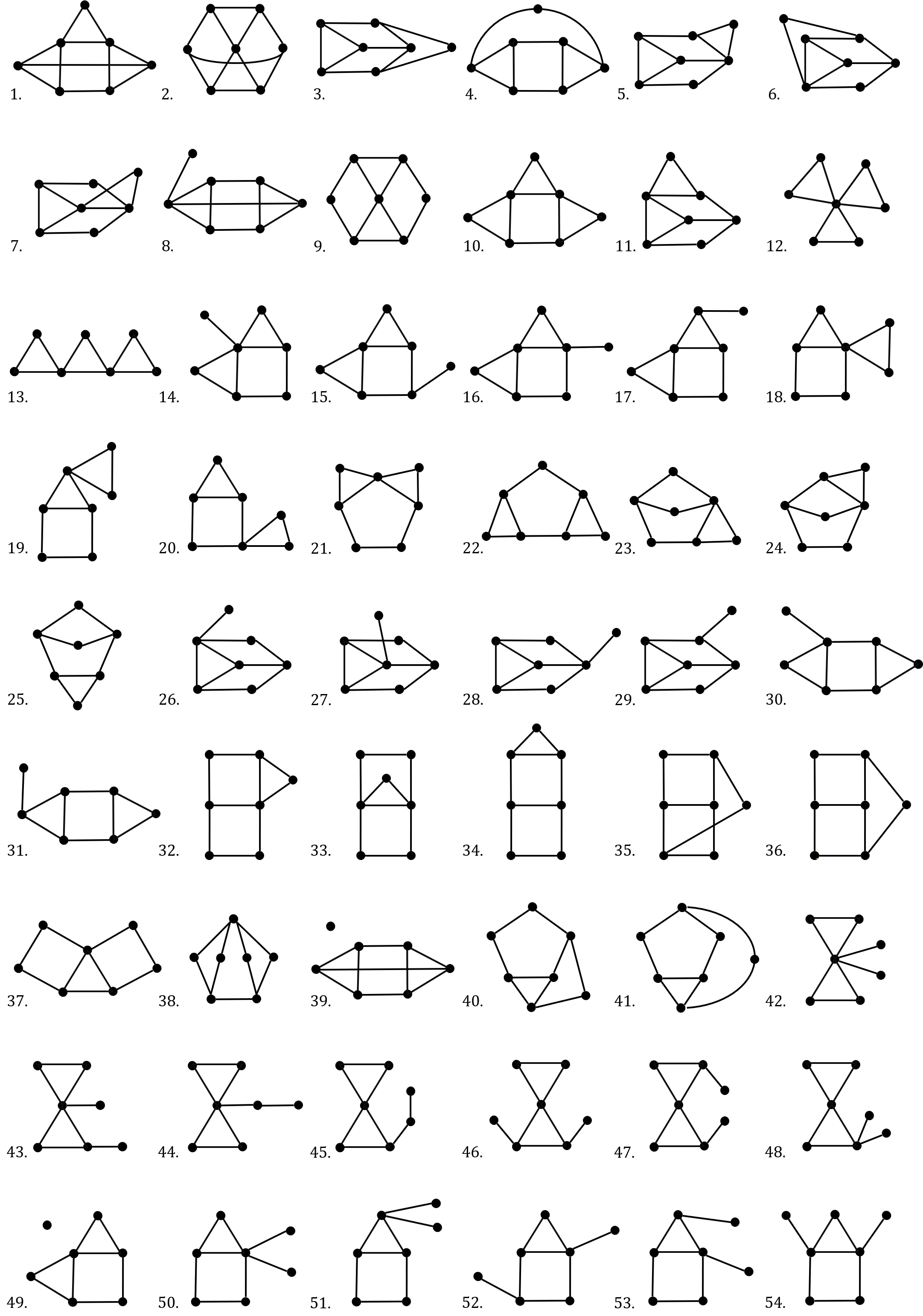}
		\centering
		\caption*{}
\end{figure}

\begin{figure}
	\includegraphics[width=1.0\textwidth]{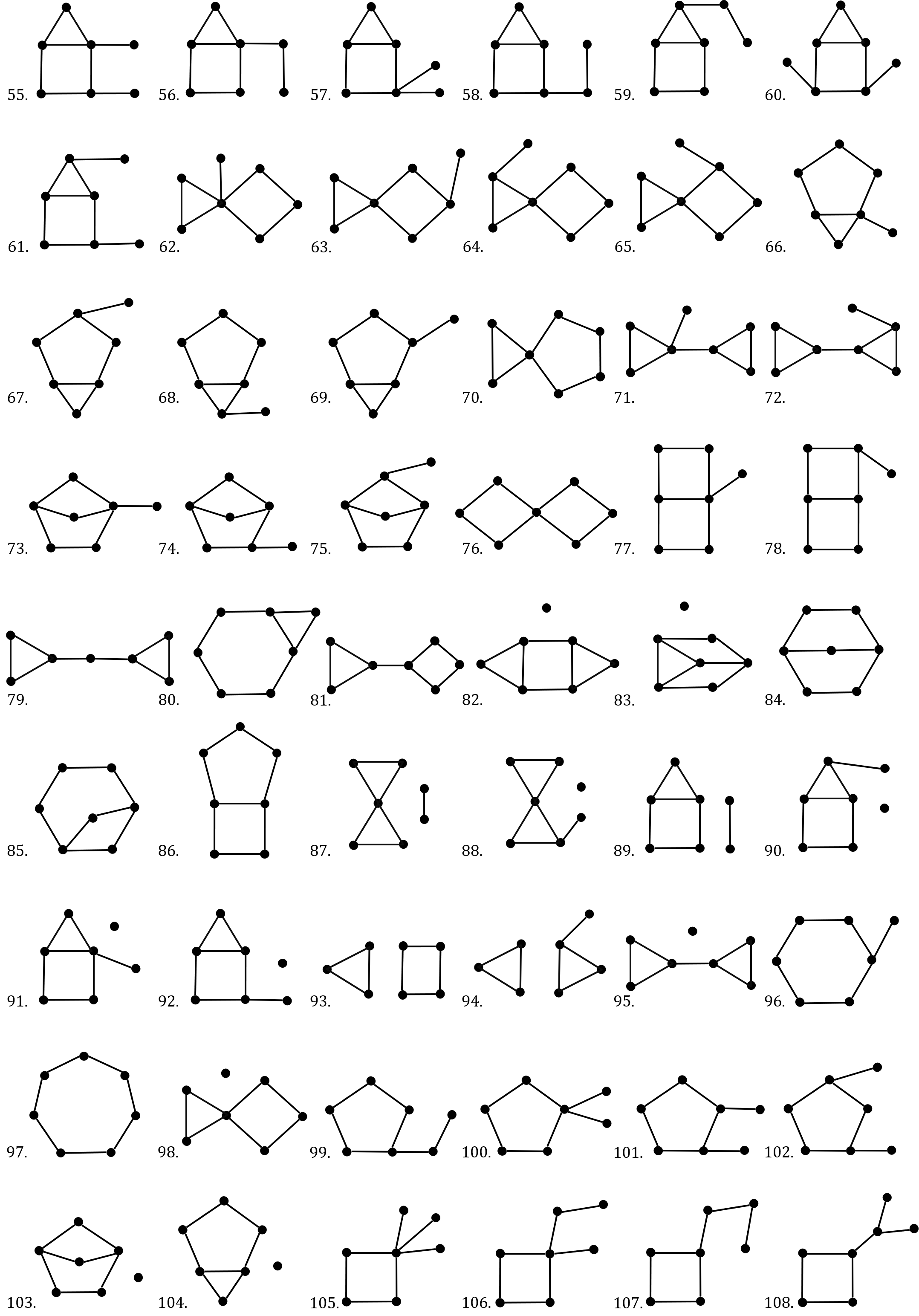}
		\centering
		\caption*{}
\end{figure}

\begin{figure}
	\includegraphics[width=1.0\textwidth]{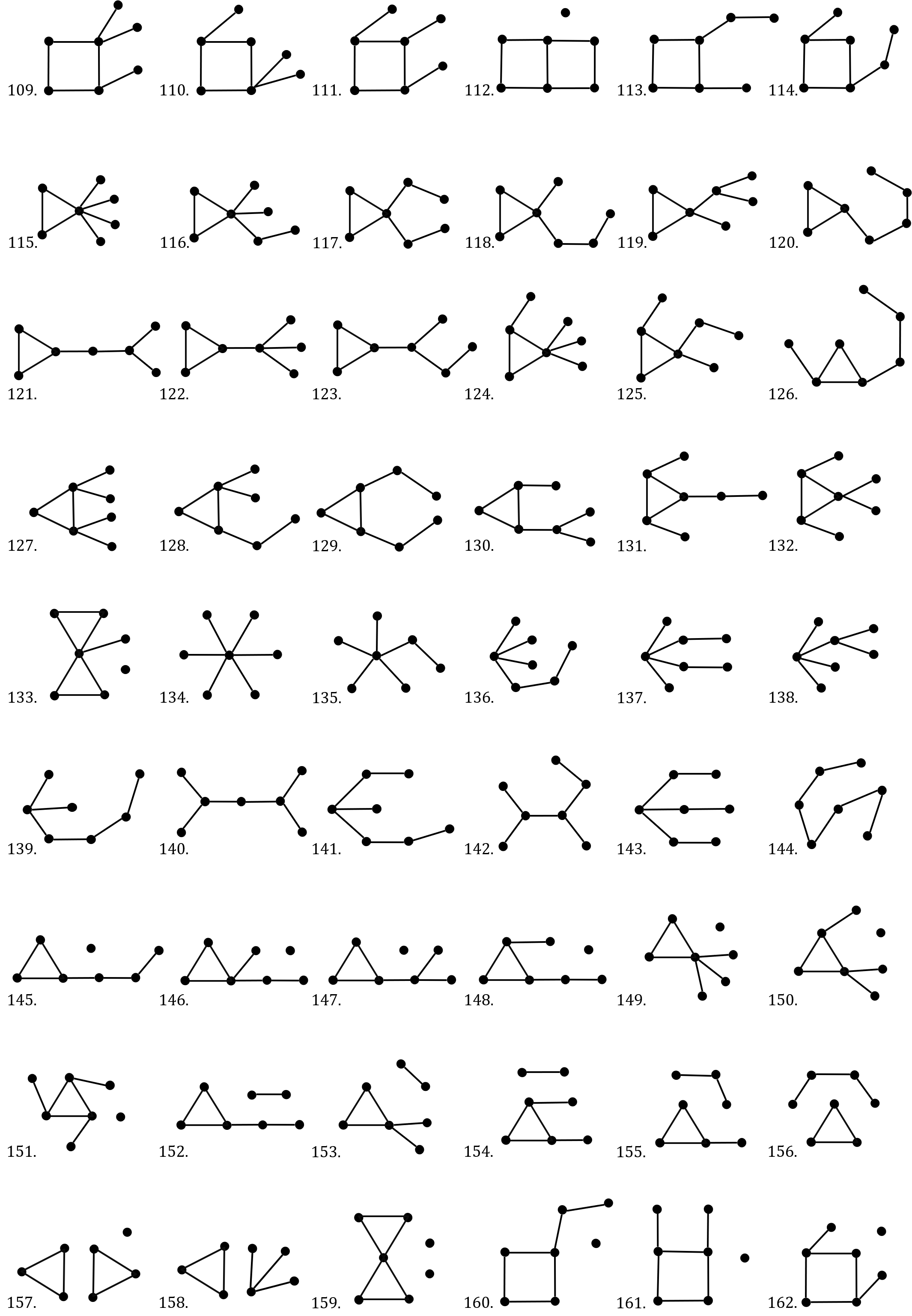}
		\centering
		\caption*{}
\end{figure}

\begingroup
\allowdisplaybreaks
 \begin{flalign*}
 z_1=&\frac{1}{4}nk(k-2)-n_3,\\
 z_2=&n_3-\frac{z_{11}}{4},\\
 z_3=&z_{11}-2n_3,\\
 z_4=&2n_3,\\
 z_5=&4n_3,\\
 z_6=&\frac{z_{11}}{2},\\
 z_7=&2n_3,\\
 z_8=&\frac{1}{2}nk(k-2)(k-4)-2(k-4)n_3,\\
 z_9=&\frac{1}{4}nk(k-2)-n_3+\frac{z_{11}}{4},\\
 z_{10}=&2n_3,\\
 z_{11}=&z_{11},\\
z_{12}=&\frac{1}{48}nk(k-2)(k-4),\\
z_{13}=&\frac{1}{8}nk(k-2)(k-4)-n_3,\\
z_{14}=&\frac{1}{2}nk(k-2)(k-4),\\
z_{15}=&\frac{1}{2}nk(k-2)(k-4)-2n_3-\frac{z_{11}}{2},\\
z_{16}=&nk(k-2)(k-4)-4n_3,\\
z_{17}=&nk(k-2)(k-4)-z_{11},\\
z_{18}=&\frac{1}{2}nk(k-2)(k-4),\\
z_{19}=&\frac{1}{4}nk(k-2)(k-4)-\frac{z_{11}}{2},\\
z_{20}=&\frac{1}{2}nk(k-2)(k-4)-4n_3,\\
z_{21}=&nk(k-2)(k-3)(k-4)-2n_3-\frac{z_{11}}{2},\\
z_{22}=&nk(k-2)(k-4)-8n_3,\\
z_{23}=&nk(k-2)(k-4)-4n_3,\\
z_{24}=&2nk(k-2)(k-4)-8n_3-z_{11},\\
z_{25}=&\frac{1}{2}nk(k-2)(k-4)-z_{11},\\
z_{26}=&4(k-4)n_3-z_{11},\\
z_{27}=&2(k-4)n_3,\\
z_{28}=&2(k-6)n_3,\\
z_{29}=&4(k-4)n_3-2z_{11},\\
z_{30}=&4(k-5)n_3,\\
z_{31}=&2(k-6)n_3,\\
z_{32}=&nk(k-2)(k-4)-4n_3,\\
z_{33}=&\frac{1}{4}nk(k-2)(k-4)-n_3,\\
z_{34}=&\frac{1}{2}nk(k-2)(k-4)-\frac{3}{2}z_{11},\\
z_{35}=&nk(k-2)(k-4)-2n_3-\frac{3}{2}z_{11},\\
z_{36}=&\frac{1}{6}nk(k-2)(k-4)-\frac{z_{11}}{2},\\
z_{37}=&\frac{1}{2}nk(k-2)(k-4)-\frac{z_{11}}{2},\\
z_{38}=&\frac{1}{4}nk(k-2)(k-4)-n_3-\frac{z_{11}}{4},\\
z_{39}=&\frac{1}{24}nk(k-2)(k-4)(k-8)-\frac{1}{3}(n-6k+15)n_3,\\
z_{40}=&2nk(k-2)(k-4)-8n_3-2z_{11},\\
z_{41}=&\frac{2}{3}nk(k-2)(k-4)-4n_3,\\
z_{42}=&\frac{1}{16}nk(k-2)(k-4)(k-6),\\
z_{43}=&\frac{1}{2}nk(k-2)(k-4)(k-5),\\
z_{44}=&\frac{1}{8}nk(k-2)(k-4)(k-6),\\
z_{45}=&\frac{1}{2}nk(k-2)(k-4)(k-8)+4n_3+z_{11},\\
z_{46}=&\frac{1}{4}nk(k-2)(k-4)(k-5)+\frac{z_{11}}{2},\\
z_{47}=&\frac{1}{2}nk(k-2)(k-4)(k-6)+2n_3+\frac{z_{11}}{2},\\
z_{48}=&\frac{1}{4}nk(k-2)(k-4)(k-6)+2n_3,\\
z_{49}=&\frac{1}{4}nk(k-2)(k-4)(k-8)+2n_3+\frac{z_{11}}{2},\\
z_{50}=&\frac{1}{2}nk(k-2)(k-4)(k-6),\\
z_{51}=&\frac{1}{4}nk(k-2)(k-4)^2-2(k-4)n_3+\frac{z_{11}}{2},\\
z_{52}=&nk(k-2)(k-4)(k-6)-4(k-6)n_3+z_{11},\\
z_{53}=&nk(k-2)(k-4)^2-4(k-4)n_3+z_{11},\\
z_{54}=&\frac{1}{2}nk(k-2)(k-4)(k-5)+2n_3,\\
z_{55}=&nk(k-2)(k-4)(k-5)-4(k-5)n_3,\\
z_{56}=&nk(k-2)(k-4)(k-7)+8n_3+z_{11},\\
z_{57}=&\frac{1}{2}nk(k-2)(k-4)(k-6)-4(k-6)n_3,\\
z_{58}=&nk(k-2)(k-4)(k-8)-4(k-8)n_3+3z_{11},\\
z_{59}=&\frac{1}{2}nk(k-2)(k-4)(k-7)-2(k-6)n_3+2z_{11},\\
z_{60}=&\frac{1}{2}nk(k-2)(k-4)(k-5)-2(2k-9)n_3+\frac{3}{2}z_{11},\\
z_{61}=&nk(k-2)(k-4)(k-5)-8(k-5)n_3+2z_{11},\\
z_{62}=&\frac{1}{4}nk(k-2)(k-4)(k-6),\\
z_{63}=&\frac{1}{4}nk(k-2)(k-4)(k-8)+4n_3,\\
z_{64}=&\frac{1}{2}nk(k-2)(k-4)(k-6)+4n_3,\\
z_{65}=&\frac{1}{2}nk(k-2)(k-4)(k-6)+4n_3,\\
z_{66}=&2nk(k-2)(k-4)(k-6)-4(k-7)n_3+z_{11},\\
z_{67}=&nk(k-2)(k-4)(k-8)-2(k-14)n_3,\\
z_{68}=&nk(k-2)(k-4)(k-6)-2(k-9)n_3+z_{11},\\
z_{69}=&2nk(k-2)(k-4)(k-7)-4(k-9)n_3+4z_{11},\\
z_{70}=&\frac{1}{2}nk(k-2)(k-4)(k-8)+6n_3+\frac{z_{11}}{2},\\
z_{71}=&\frac{1}{4}nk(k-2)(k-4)(k-6)-2(k-6)n_3,\\
z_{72}=&\frac{1}{2}nk(k-2)(k-4)(k-7)-4(k-7)n_3,\\
z_{73}=&nk(k-2)(k-4)(k-7)-4(k-7)n_3+z_{11},\\
z_{74}=&nk(k-2)(k-4)(k-7)-4(k-6)n_3+3z_{11},\\
z_{75}=&nk(k-2)(k-4)(k-6)-2(2k-11)n_3+2z_{11},\\
z_{76}=&\frac{1}{8}nk(k-2)(k-4)(k-8)+n_3+\frac{z_{11}}{4},\\
z_{77}=&\frac{1}{2}nk(k-2)(k-4)(k-6)-2(k-6)n_3,\\
z_{78}=&nk(k-2)(k-4)(k-7)-4(k-6)n_3+4z_{11},\\
z_{79}=&\frac{1}{8}nk(k-2)(k-4)(k-12)+6n_3,\\
z_{80}=&\frac{1}{2}nk(k-2)(k-4)(2k-7)+16n_3+\frac{3}{2}z_{11},\\
z_{81}=&\frac{1}{4}nk(k-2)(k-4)(k-10)+2n_3+2z_{11},\\
z_{82}=&(n-6k+18)n_3,\\
z_{83}=&2(n-6k+16)n_3+z_{11},\\
z_{84}=&\frac{1}{2}nk(k-2)(k-4)(2k-7)+16n_3+\frac{3}{2}z_{11},\\
z_{85}=&\frac{1}{2}nk(k-2)(k-4)(k-9)+8n_3+\frac{3}{2}z_{11},\\
z_{86}=&nk(k-2)(k-4)(k-8)+12n_3+\frac{5}{2}z_{11},\\
z_{87}=&\frac{1}{32}nk(k-2)(k-4)(k^2-16k+68)-n_3-\frac{z_{11}}{4},\\
z_{88}=&\frac{1}{4}nk(k-2)(k-4)(k^2-12k+40)-4n_3-z_{11},\\
z_{89}=&\frac{1}{8}nk(k-2)(k-4)(k^2-16k+68)+2(k-8)n_3-\frac{3}{2}z_{11},\\
z_{90}=&\frac{1}{4}nk(k-2)(k-4)(k^2-11k+32)-2(n-7k+21)n_3-2z_{11},\\
z_{91}=&\frac{1}{2}nk(k-2)(k-4)(k^2-12k+38)+4(k-6)n_3-z_{11},\\
z_{92}=&\frac{1}{2}nk(k-2)(k-4)(k^2-12k+40)-4(n-7k+24)n_3-3z_{11},\\
z_{93}=&\frac{1}{96}nk(k-2)(k-4)(k^2-20k+114)-n_3-\frac{z_{11}}{2},\\
z_{94}=&\frac{1}{24}nk(k-2)(k-4)(k-8)(k-12)+2(k-8)n_3,\\
z_{95}=&\frac{1}{16}nk(k-2)(k-4)(k^2-12k+44)-(n-6k+21)n_3,\\
z_{96}=&\frac{1}{2}nk(k-2)(k-4)(2k^2-29k+118)+(6k-80)n_3-\frac{11}{2}z_{11},\\
z_{97}=&\frac{1}{14}nk(k-2)(k-4)(2k^2-30k+133)-10n_3-z_{11},\\
z_{98}=&\frac{1}{4}nk(k-2)(k-4)(n-6k+20)-4n_3,\\
z_{99}=&nk(k-2)(k-4)(k^2-14k+56)+(6k-76)n_3-6z_{11},\\
z_{100}=&\frac{1}{2}nk(k-2)(k-4)(k^2-14k+50)+(6k-40)n_3-z_{11},\\
z_{101}=&nk(k-2)(k-4)(k^2-13k+44)+12(k-6)n_3-5z_{11},\\
z_{102}=&nk(k-2)(k-4)(k^2-14k+53)+(12k-94)n_3-\frac{9}{2}z_{11},\\
z_{103}=&\frac{1}{2}nk(k-2)(k-4)(n-6k+21)-2(n-6k+20)n_3-2z_{11},\\
z_{104}=&nk(k-2)(k-4)(n-6k+21)-2(n-6k+26)n_3-2z_{11},\\
z_{105}=&\frac{1}{12}nk(k-2)(k-4)(k-6)(k-8),\\
z_{106}=&\frac{1}{2}nk(k-2)(k-4)(k^2-13k+44)+4(k-7)n_3-z_{11},\\
z_{107}=&\frac{1}{2}nk(k-2)(k-4)(k^2-15k+64)+4(k-11)n_3-6z_{11},\\
z_{108}=&\frac{1}{4}nk(k-2)(k-4)(k^2-14k+52)+4(k-6)n_3-\frac{5}{2}z_{11},\\
z_{109}=&\frac{1}{2}nk(k-2)(k-4)(k-6)^2+4(k-6)n_3,\\
z_{110}=&\frac{1}{4}nk(k-2)(k-4)(k^2-14k+50)+2(2k-13)n_3-\frac{z_{11}}{2},\\
z_{111}=&\frac{1}{2}nk(k-2)(k-4)(k^2-12k+37)+8(k-5)n_3-\frac{5}{2}z_{11},\\
z_{112}=&\frac{1}{8}nk(k-2)(k-4)(k^2-12k+44)-(n-6k+20)n_3-\frac{5}{4}z_{11},\\
z_{113}=&nk(k-2)(k-4)(k^2-13k+46)+4(3k-20)n_3-6z_{11},\\
z_{114}=&\frac{1}{2}nk(k-2)(k-4)(k^2-14k+54)+8(k-7)n_3-4z_{11},\\
z_{115}=&\frac{1}{48}nk(k-2)(k-4)(k-6)(k-8),\\
z_{116}=&\frac{1}{4}nk(k-2)(k-4)(k-6)^2,\\
z_{117}=&\frac{1}{4}nk(k-2)(k-4)(k^2-12k+41)-6n_3-\frac{z_{11}}{2},\\
z_{118}=&\frac{1}{2}nk(k-2)(k-4)(k^2-14k+51)+4(k-8)n_3-z_{11},\\
z_{119}=&\frac{1}{4}nk(k-2)(k-4)(k^2-12k+36)+2(k-6)n_3,\\
z_{120}=&\frac{1}{2}nk(k-2)(k-4)(k^2-16k+71)+4(k-14)n_3-4z_{11},\\
z_{121}=&\frac{1}{4}nk(k-2)(k-4)(k^2-16k+68)+4(k-10)n_3,\\
z_{122}=&\frac{1}{12}nk(k-2)(k-4)(k-6)(k-8)+2(k-6)n_3,\\
z_{123}=&\frac{1}{2}nk(k-2)(k-4)(k^2-14k+52)+4(2k-13)n_3-4z_{11},\\
z_{124}=&\frac{1}{6}nk(k-2)(k-4)(k-5)(k-6),\\
z_{125}=&nk(k-2)(k-4)(k^2-11k+32)+4(k-7)n_3-z_{11},\\
z_{126}=&nk(k-2)(k-4)(k^2-14k+52)+12(k-8)n_3-3z_{11},\\
z_{127}=&\frac{1}{8}nk(k-2)(k-4)(k^2-10k+26)-n_3,\\
z_{128}=&\frac{1}{2}nk(k-2)(k-4)(k^2-12k+38)+4(k-6)n_3-z_{11},\\
z_{129}=&\frac{1}{2}nk(k-2)(k-4)(k^2-13k+49)+2(2k-17)n_3-\frac{7}{2}z_{11},\\
z_{130}=&\frac{1}{2}nk(k-2)(k-4)(k^2-13k+44)+8(k-6)n_3-2z_{11},\\
z_{131}=&\frac{1}{2}nk(k-2)(k-4)(k^2-11k+32)+2(3k-16)n_3-2z_{11},\\
z_{132}=&\frac{1}{4}nk(k-2)(k-4)^2(k-5)+2(k-4)n_3-\frac{z_{11}}{2},\\
z_{133}=&\frac{1}{8}nk(k-2)(k-4)(n-6k+17),\\
z_{134}=&\frac{1}{720}nk(k-2)(k-4)(k-6)(k-8)(k-10),\\
z_{135}=&\frac{1}{24}nk(k-2)(k-4)(k-6)^2(k-8),\\
z_{136}=&\frac{1}{6}nk(k-2)(k-4)(k^3-20k^2+141k-348)-(8k-52)n_3+z_{11},\\
z_{137}=&\frac{1}{4}nk(k-2)(k-4)(k^3-18k^2+113k-250)-2(3k-20)n_3+z_{11},\\
z_{138}=&\frac{1}{12}nk(k-2)(k-4)(k-6)(k^2-12k+38)-2(k-6)n_3,\\
z_{139}=&\frac{1}{2}nk(k-2)(k-4)(k^3-20k^2+147k-396)-4(7k-55)n_3+9z_{11},\\
z_{140}=&\frac{1}{8}nk(k-2)(k-4)(k^3-20k^2+144k-372)-8(k-7)n_3+\frac{3}{2}z_{11},\\
z_{141}=&nk(k-2)(k-4)(k^3 - 19k^2 + 133k-344)-4(12k-89)n_3 + 20z_{11},\\
z_{142}=&\frac{1}{2}nk(k-2)(k-4)(k^3-18k^2+116k-266)- 8(3k -17)n_3+7z_{11},\\
z_{143}=&\frac{1}{6}nk(k-2)(k-4)(k^3-18k^2+120k-302)-6(k-8)n_3+3z_{11},\\
z_{144}=&\frac{1}{2}nk(k-2)(k-4)(k^3-20k^2+150k-426)-10(2k-21)n_3+\frac{27}{2}z_{11},\\
z_{145}=&\frac{1}{4}nk(k-2)(k-4)(k-8)(k^2-11k+44)+4n_3(n-7k+34)+4z_{11},\\
z_{146}=&\frac{1}{4}nk(k-2)(k-4)(k^3-17k^2+104k-232)-4(k-8)n_3+z_{11},\\
z_{147}=&\frac{1}{8}nk(k-2)(k-4)(k^3 -18k^2+118k -284)+2n_3(n-8k+32)+2z_{11},\\
z_{148}=&\frac{1}{2}nk(k-2)(k-4)(k^3 -17k^2 +106k - 248)+4n_3(n-9k+40)+7z_{11},\\
z_{149}=&\frac{1}{12}nk(k-2)(k-4)(k-6)(n-6k+18),\\
z_{150}=&\frac{1}{4}nk(k-2)(k-4)(k - 6)(k^2-10k+30)-4n_3(k-6)+z_{11},\\
z_{151}=&\frac{1}{12}nk(k-2)(k-4)(k^3-15k^2+80k-152)+\frac{2}{3}n_3(n-9k+30)+z_{11},\\
z_{152}=&\frac{1}{8}nk(k-2)(k-4)(k^3-20k^2+148k-412)-2n_3(3k-26)+4z_{11},\\
z_{153}=&\frac{1}{16}nk(k-2)(k-4)(k^3-20k^2+140k-344)-2(k-7)n_3+\frac{z_{11}}{2},\\
z_{154}=&\frac{1}{8}nk(k-2)(k-4)(k-7)(k^2-12k+44)-2n_3(3k-20)+\frac{3}{2}z_{11},\\
z_{155}=&\frac{1}{8}nk(k-2)(k-4)(k-8)(k^2-14k+60)-2n_3(5k-39)+2z_{11},\\
z_{156}=&\frac{1}{24}nk(k-2)(k-4)(k^3-23k^2+192k -588)-2n_3(2k-17)+2z_{11},\\
z_{157}=&\frac{1}{288}nk(k-2)(k-4)(k-12)(k^2-12k+50)+\frac{1}{3}n_3(n-6k+24),\\
z_{158}=&\frac{1}{72}nk(k-2)(k-4)(k-8)(k^2-16k+78)-2(k-8)n_3,\\
z_{159}=&\frac{1}{64}nk(k-2)(k-4)(k^3 -18k^2+114k-268)+n_3+\frac{z_{11}}{4},\\
z_{160}=&\frac{1}{4}nk(k-2)(k-4)(k^3-18k^2+120k-308)+4n_3(n-7k+29)+7z_{11},\\
z_{161}=&\frac{1}{4}nk(k-2)(k-4)(k^3-17k^2+106k-242)+2n_3(n-10k+42)+\frac{7}{2}z_{11},\\
z_{162}=&\frac{1}{8}nk(k-2)(k-4)(k^3-18k^2+118k-288)+2n_3(n-8k+32)+\frac{5}{2}z_{11},\\
z_{163}=&\frac{1}{8}nk(k-2)(k-4)(k^3-18k^2+116k-268)-2n_3(2k-13)+\frac{z_{11}}{2},\\
z_{164}=&\frac{1}{8}nk(k-2)(k-4)(k^3-20k^2+148k-404)-2n_3(4k-29)+4z_{11},\\
z_{165}=&\frac{1}{32}nk(k-2)(k-4)(k^3-22k^2+ 178k - 540)-2n_3(k-9)+2z_{11},\\
z_{166}=&\frac{1}{16}nk(k-2)(k-4)(k^3 - 18 k^2+ 114k-268)+2n_3(n-6k+19)+\frac{3}{2}z_{11},\\
z_{167}=&\frac{1}{2}nk(k-2)(k-4)(k^3-18k^2+120k-304)+2n_3(3n-23k+104)+9z_{11},\\
z_{168}=&\frac{1}{20}nk(k-2)(k-4)(k^3-20k^2+152k-440)-2n_3(k-11)+\frac{3}{2}z_{11},\\
z_{169}=&\frac{1}{24}nk(k-2)(k-4)(2k^3-37k^2+257k-710)+n_3(n-6k+35)+\frac{7}{4}z_{11},\\
z_{170}=&\frac{1}{240}nk(k-2)(k-4)(k-6)(k-8)(k^2-12k+40),\\
z_{171}=&\frac{1}{12}nk(k-2)(k-4)(k^4-23k^3+208k^2-890k+1536)+4n_3(2k-13)-z_{11},\\
z_{172}=&\frac{1}{96}nk(k-2)(k-4)(k^4-26k^3+268k^2-1304k+2520)+(2k-13)n_3-\frac{z_{11}}{4},\\
z_{173}=&\frac{1}{4}nk(k-2)(k-4)(k^4-23k^3+216k^2-998k+1944)-4n_3(2n-21k+105)-15z_{11},\\
z_{174}=&\frac{1}{4}nk(k-2)(k-4)(k^4-22k^3+198k^2-878k+1644)-6n_3(n-11k+54)-\frac{29}{2}z_{11},\\
z_{175}=&\frac{1}{24}nk(k-2)(k-4)(k^4-26k^3+274k^2-1400k+2928)+2n_3(7k-50)-\frac{5}{2}z_{11},\\
z_{176}=&\frac{1}{84}nk(k-2)(k-4)(k^4-24k^3+236k^2-1136k+2264)+6(5k-34)n_3- 9z_{11},\\
z_{177}=&\frac{1}{16}nk(k-2)(k-4)(k^4-22k^3+194k^2-824k+1428)-n_3(n-14k+64)-\frac{9}{4}z_{11},\\
z_{178}=&\frac{1}{4}nk(k-2)(k-4)(k^4-23k^3+220k^2-1052k+2174)-2n_3(5n-44k+240)-21z_{11},\\
z_{179}=&\frac{1}{8}nk(k-2)(k-4)(k^4-24k^3+241k^2-1212k+2596)+6(5k-41)n_3-\frac{27}{2}z_{11},\\
z_{180}=&\frac{1}{8}nk(k-2)(k-4)(k^4-25k^3+258k^2-1316k+2820)+4(9k-70)n_3-\frac{27}{2}z_{11},\\
z_{181}=&\frac{1}{16}nk(k-2)(k-4)(k^4-22k^3+198k^2-876k+1664)-2n_3(2n-13k+52)-4z_{11},\\
z_{182}=&\frac{1}{32}nk(k-2)(k-4)(k^4 - 22k^3 + 190k^2 - 780k + 1312)+2(k-7)n_3-\frac{z_{11}}{2},\\
z_{183}=&\frac{1}{16}nk(k-2)(k-4)(k^4-21k^3+176k^2-714k+1212)-2n_3(n-8k+31)-\frac{5}{2}z_{11},\\
z_{184}=&\frac{1}{16}nk(k-2)(k-4)(k^4-24k^3+232k^2-1096k+2184)-2n_3(n-11k+59)-4z_{11},\\
z_{185}=&\frac{1}{192}nk(k-2)(k-4)(k^4-26k^3+280k^2-1508k+3456)+2n_3(k-8)-z_{11},\\
z_{186}=&\frac{1}{48}nk(k-2)(k-4)(k^4-26k^3+274k^2-1412k+3096)-2n_3(n-8k+40)-2z_{11},\\
z_{187}=&\frac{1}{16}nk(k-2)(k-4)(k^4-22k^3+198k^2-876k+1648)-4n_3(n-7k+27)-\frac{9}{2}z_{11},\\
z_{188}=&\frac{1}{64}nk(k-2)(k-4)(k^4-24k^3+238k^2-1176k+2520)-n_3(n-8k+39)-\frac{5}{2}z_{11},\\
z_{189}=&\frac{1}{40}nk(k-2)(k-4)(k^4-22k^3+202k^2-924k+1840)-2n_3(n-6k+24)-2z_{11},\\
z_{190}=&\frac{1}{192}nk(k-2)(k-4)(k^5-28k^4+326k^3-2000k^2+6576k-9408)-(2k-13)n_3+\frac{z_{11}}{4},\\
z_{191}=&\frac{1}{16}nk(k-2)(k-4)(k^5-26k^4+290k^3-1764k^2+5952k-9072)\\&+2n_3(4n-37k+166)+12z_{11},\\
z_{192}=&\frac{1}{48}nk(k-2)(k-4)(k^5-28k^4+338k^3-2228k^2+8136k-13320)\\&+(k^2-30k+166)n_3+5z_{11},\\
z_{193}=&\frac{1}{16}nk(k-2)(k-4)(k^5-26k^4+295k^3-1858k^2+6594k-10812)\\&+2n_3(7n-52k+236)+\frac{37}{2}z_{11},\\
z_{194}=&\frac{1}{16}nk(k-2)(k-4)(k^5-27k^4+320k^3-2108k^2+7824k-13312)\\&+2n_3(4n-48k+269)+23z_{11},\\
z_{195}=&\frac{1}{64}nk(k-2)(k-4)(k^5-28k^4+348k^3-2420k^2+9496k-16912)\\&-2n_3(11k-81)+\frac{15}{2}z_{11},\\
z_{196}=&\frac{1}{96}nk(k-2)(k-4)(k^5-26k^4+286k^3-1700k^2+5592k-8400)\\&+2n_3(n-7k+28)+2z_{11},\\
z_{197}=&\frac{1}{192}nk(k-2)(k-4)(k^5-28k^4+338k^3-2240k^2+8304k-14208)\\&+2n_3(n-7k+31)+2z_{11},\\
z_{198}=&\frac{1}{384}nk(k-2)(k-4)(k^5-26k^4+292k^3-1808k^2+6276k-10104)\\&+n_3(n-6k+22)+z_{11},\\
z_{199}=&\frac{1}{2304}nk(k-2)(k-4)(k^6-30k^5+396k^4-3000k^3+14084k^2-39768k+54528)\\&-\frac{1}{3}n_3(2n-12k+45)-\frac{z_{11}}{2},\\
z_{200}=&\frac{1}{288}nk(k-2)(k-4)(k^6-30k^5+396k^4-3000k^3+14012k^2-38952k+51312)\\&-2n_3(n-9k+40)-\frac{5}{2}z_{11},\\
z_{201}=&\frac{1}{96}nk(k-2)(k-4)(k^6-29k^5+376k^4-2846k^3+13500k^2-38664k+53472)\\&-2n_3(5n-36k+153)-\frac{23}{2}z_{11},\\
z_{202}=&\frac{1}{64}nk(k-2)(k-4)(k^6-30k^5+406k^4-3228k^3+16128k^2-48672k+70624)\\&-2n_3(7n-63k+306)-22z_{11},\\
z_{203}=&\frac{1}{384}nk(k-2)(k-4)(k^6-30k^5+412k^4-3364k^3+17392k^2-54576k+82272)\\&-\frac{1}{3}n_3(4n-60k+345)-\frac{19}{4}z_{11},\\
z_{204}=&\frac{1}{768}nk(k-2)(k-4)(k^7-32k^6+464k^5-4048k^4+23428k^3-91424k^2\\&+226032k-279168)+2n_3(3n-21k+86)+6z_{11},\\
z_{205}=&\frac{1}{768}nk(k-2)(k-4)(k^7-32k^6+470k^5-4204k^4+25168k^3-102104k^2\\&+263184k-338976)+n_3(6n-48k+215)+\frac{31}{4}z_{11},\\
z_{206}=&\frac{1}{7680}nk(k-2)(k-4)(k^8-34k^7+528k^6-5040k^5+33236k^4-157832k^3\\&+532944k^2-1184160k+1353600)-n_3(3n-20k+78)-\frac{11}{4}z_{11},\\
z_{207}=&\frac{1}{322560}nk(k-2)(k-4)(k^9-36k^8+586k^7-5856k^6+41156k^5-217104k^4\\&+869528k^3-2564928k^2+5119200k-5376000)+\frac{n_3}{3}(n-6k+21)+\frac{z_{11}}{4},\\
z_{208}=&\frac{1}{32}nk(k-2)(k-4)(k^5-28k^4+342k^3-2308k^2+8724k-15016)\\&+n_3(5n-58k+321)+\frac{45}{4}z_{11}.\\
\end{flalign*}
\endgroup

\section{Conclusion}

In conclusion, we have studied the structure of a class of strongly regular graphs with parameters $\lambda = 1$ and $\mu = 2$. In particular, we have determined all induced subgraphs of order seven and computed their relative frequencies. In doing so, we have used extensively Wolfram Alpha, a computational engine, for algebraic calculations \cite{Wolfram}. These results provide a more detailed understanding of the local structure of such graphs and illustrate how global parameters influence small-scale configurations.

Although the existence of several graphs in this parameter regime remains unresolved—most notably the $srg(99,14,1,2)$, known as the 99-vertex Conway graph problem \cite{Conway}—we believe that the results developed in this paper may offer useful insights toward addressing such open questions. In particular, the classification and frequency analysis of small induced subgraphs could serve as a tool for constraining or guiding future constructions.


\end{document}